\documentclass{amsart}

\usepackage{amssymb, amsmath, amsthm}
\usepackage{hyperref}

\title{Modular equations of order $p$ and theta functions}

\author[Kopeliovich]{Yaacov Kopeliovich}



\newtheorem{thm}{Theorem}[section]
\newtheorem{cor}[thm]{Corollary}
\newtheorem{lem}[thm]{Lemma}

\theoremstyle{definition}
\newtheorem{defn}[thm]{Definition}
\theoremstyle{Assumption}

\theoremstyle{remark}

\begin{document}
\begin{abstract}
Let $p$ be a prime integer and $\mathbf{H_g}$ be a collection of complex  positive definite symmetric $g\times g$ matrices $\tau$. Denote by  $p\tau$ the multiplication of $\tau$ by $p.$ In this note we describe an explicit process to obtain algebraic identities  between theta functions with integral characteristics evaluated at $\tau$ and $p\tau.$
For $g=1$ this produces  modular equations between  $\lambda(\tau),\lambda(p\tau)$ where $\lambda(\tau)$ is the invariant associated with elliptic curve generated by $\tau,$ described by the equation: $y^2=x(x-1)\left(x-\lambda(\tau_1)\right).$ Consequently, if $g>1$  The algebraic identities we obtain might serve as a higher dimensional generalization for the one dimensional modular equations.
\end{abstract}

\maketitle

\medskip

\begin{center}
\emph{Dedicated to Mike Fried on his 65-th birthday and constant mathematical inspiration.}
\end{center}

\bigskip

\section{Introduction}
Let $\tau_1$ be a complex number such that $Im(\tau_1)>0$ and $Z_{\tau_1}$ is the lattice generated  by $\left\{1,\tau_1\right\}$ in $\mathbb{C}$. Let $C_1=\mathbb{C}/Z_{\tau_1},$ be the corresponding analytic elliptic curve. The algebraic equation of this curve is:
\begin{equation}
y^2=x(x-1)\left(x-\lambda(\tau_1)\right).
\end{equation}
$\lambda(\tau_1)$ is the invariant corresponding to $\tau_1$ in this equation.

\begin{defn}Let  $p$ be a prime number. A modular equation of order $p$ for $\lambda,$ is an algebraic equation between $\lambda(\tau_1)$ and $\lambda(p\tau_1).$
\end{defn}

These equations appeared naturally in the theory of elliptic integrals because, such an equation describes algebraically the analytical multiplication by $p$ on the lattice $Z_{\tau_1}$. In equivalent more contemporary terms this equation describes $\lambda$ invariant of curves $C_2$ and $\phi:C_1\mapsto C_2$ is an isogeny (finite homomorphism) of order $p.$ Equations of this type have an important role in Galois theory of complex multiplications since, if $C_1$ is a curve with complex multiplication then $\lambda(p\tau_1)$ generates interesting field extension of $Q^{ab}(\tau_1).$ Another application of one dimensional modular equations is algorithms for rapid calculation of $\pi,$ \cite{Bo}.

In recent years there are applications of modular equations to point counting algorithms of elliptic curve above finite fields $F_{p^i}.$ Modular equations of order $p$ are used to compute explicit canonical lifting of elliptic
curves above finite fields $F_{p^i},$ of characteristics $p$ to the
corresponding $p$ - adic field above it. Using the lifting we calculate the  trace of the Frobenius operator on the  $p$- adic field. Applying fixed point formulas, we find  the number of  points of elliptic curves above finite fields $F_{p^i}.$ \cite{Ma} explains the general framework for this type of algorithms.

For $p=2$ the modular equation of level $2$ is the arithmetic geometric mean (AGM) in disguise. If $a_0=a, b_0=b$ are real numbers we define the AGM iteration as:

\begin{equation}
a_n = \frac{a_n+b_n}{2}
\end{equation}
\begin{equation}
b_n=\sqrt{a_nb_n}
\end{equation}

This iteration has a strong link  with  the modular equations of order $2.$ These sequences converge to a common limit denoted by $AGM(a,b)$ see \cite{Bo}. Mestre \cite{Me}  used the AGM iteration to suggest an algorithm for point counting  defined over $F_{2^i}.$ \cite{Me} Mestre uses a higher dimensional analogue of the AGM and produces an algorithm to count the number of points of hyper-elliptic curves above fields of characteristics $2$. In a recent work Lubicz Carls and Kohel \cite{CKL} generalized Mestre's  method further in a different direction. They construct curves that according to them have good cryptographic properties using an iteration of order $3$. In one dimensional case their iteration seems to be closely related to modular equations of order $3.$ These results motivate the question whether there exists a theory of modular equations for $\tau$ an element of $\mathbf{H_g}.$ To suggest a possible generalization  recall that $\lambda(\tau_1)$ is a quotient of theta functions that is:
$$\lambda(\tau_1) =\frac{\Theta^4\left[\begin{array}{c}
         0 \\
         1
\end{array}\right]\left(0,\tau_1\right)}{\Theta^4\left[\begin{array}{c}
         0 \\
         0
\end{array}\right]\left(0,\tau_1\right)}.
$$

These are analytic functions that we define in the first section of the paper. $\lambda(p\tau_1)$ is the same quotient evaluated at $p\tau_1.$ Hence modular equations become identities between theta functions evaluated at $\tau_1$ and those evaluated at $p\tau_1.$ While an elementary dimension argument shows that analogues of $\lambda(\tau_1)$ do not exist for general $\tau \in \mathbf{H_g}$ theta functions do. Thus the question of modular equations is reduced for finding a way to produce certain identities between higher dimensional theta functions evaluated at $\tau$ and $p\tau.$ In this note we apply methods from [Ko] to suggest such a procedure. The procedure constructs equations between

\[\Theta\left[\begin{array}{c}
         \eta_i \\
         \epsilon_i
       \end{array}\right]\left(0,\tau\right)
\textit{ and  }
\Theta\left[\begin{array}{c}
         \eta_i \\
         \epsilon_i
\end{array}\right]\left(0,p\tau\right),
 \]
for any $p>2$ and $\eta_i, \epsilon_i$ are $g$ integral characteristics. We stress that in addition to applications similar in the one dimensional case and possibly in cryptography we believe that an existence of such a procedure should lead to other applications that are not present in the 1-dimensional case. For example such a procedure should produce algebraic conditions that characterize Abelian varieties that are isogenous to multiplication of elliptic curves. The problem treated in this note is addressed in the literature with  different approaches. Modular equations for hyper elliptic curves through $l$ torsion subgroups of Jacobians were defined and treated in \cite {GS}. \cite{CKL} treated the  case for $p=3$ and \cite {CL} seem to treat the general $p$ along the same lines. Their work produces identities that use the theory of algebraic theta functions and Riemann's theta formula. The relation of these identities to the current work is not obvious and we plan to investigate it further in the future.

We review the structure of this note: The first section explains the process to obtain identities between theta functions with integral characteristics at $\tau$ and $p\tau$. We apply these results in second section to the one dimensional case and explain how this leads to a proof of existence of modular equations to the $\lambda$ function. This serves as an alternative to the classical theory of modular polynomials that has no analogue in the higher dimensional case. The last section we produce modular equations
for $p=3,7.$

\bigskip

\noindent \textbf{Acknowledgements: } We thank David Lubicz whose interest in these questions prompted the initial motivation for this work. We thank David Kohel and Hershel Farkas for reading and providing valuable suggestions on an earlier version for this note. We especially thank the referee for very constructive remarks that substantially improved the  presentation of this note. This work was partially done while the author visited Oakland university and the author thanks the department for support and hospitality.

\section{Modular equations}
We remind the reader the definition and main properties of theta functions:

\begin{defn}
Let $\tau$ be a complex $g\times g$ matrix such that:

\begin{itemize}
\item $\tau=\tau^t$ i.e. $\tau$ is symmetric
\item $Im\tau$ is a positive definite quadratic form.
\end{itemize}

Then $\tau \in \mathbf{H_g}.$ Let $\left[\begin{array}{c}
         \varepsilon \\
         \varepsilon' \\
       \end{array}\right]$ be a real $2g$ vector. \textbf{Theta function} is a complex analytic function on $C^g\times
\mathbf{H_g}$ such that:

\[\Theta\left[
\begin{array}{c}
  \varepsilon \\
  \varepsilon' \\
\end{array}
\right](z, \tau) = \sum_{l\in \mathbb{Z}^{g}} exp 2\pi i
\left\{\frac{1}{2}\left(l+\frac{\varepsilon}{2}\right)^t\tau\left(l+\frac{\varepsilon}{2}\right)+
\left(l+\frac{\epsilon}{2}\right)^t\left(z+\frac{\epsilon'}{2}\right)\right\}
\]
\end{defn}

Note that the definition is the classical definition of theta function. The modern authors omit the factor $\frac{1}{2}.$ We list the main properties of theta functions:\\

\begin{itemize}
  \item $\Theta\left[
\begin{array}{c}
  \varepsilon + 2m\\
  \varepsilon'+ 2e \\
\end{array}\right] (z, \tau)= exp\pi i\left\{\varepsilon^t e\right\}\Theta\left[
\begin{array}{c}
  \varepsilon \\
  \varepsilon' \\
\end{array}\right](z, \tau)$ and $m,e \in Z^g$

  \item $\Theta\left[
\begin{array}{c}
  \varepsilon \\
  -\varepsilon' \\
\end{array}\right](z, \tau) = \Theta\left[
\begin{array}{c}
  \varepsilon \\
  \varepsilon' \\
\end{array}\right](-z, \tau)$

\item $\Theta\left[\begin{array}{c}
  \varepsilon \\
  \varepsilon' \\
\end{array}\right](z + n + m^t\tau, \tau) =  \exp {2\pi i\left\{{\frac{n^t
\varepsilon-m^t \varepsilon'}{2}-m^tz-m^t\tau
m}\right\}}\Theta\left[\begin{array}{c}
  \varepsilon \\
  \varepsilon' \\
\end{array}\right](z,\tau)$\\
\end{itemize}

For the proof of these properties that follow from a careful series manipulation see \cite{Mu}, \cite{RF} or \cite{Ko}. We remind the reader the notion of integral (rational) theta characteristics:
\begin{defn} The functions
$$\Theta\left[
\begin{array}{c}
  \varepsilon \\
  \varepsilon' \\
\end{array}
\right](z, \tau)$$ are called theta functions with \textbf{integral} (\textbf{rational})  characteristics if $\left[
\begin{array}{c}
  \varepsilon \\
  \varepsilon' \\
\end{array}
\right] \in \mathbb{Z}^{2g}\left(\mathbb{Q}^{2g}\right)$
\end{defn}

From property (2) we see that we can assume that
$\left[\begin{array}{c}
  \varepsilon \\
  \varepsilon' \\
\end{array}\right] \in Z_2^{2g}$
Further $\Theta\left[
\begin{array}{c}
  \varepsilon \\
  \varepsilon' \\
\end{array}
\right](z, \tau)$ will be even or odd if  the scalar product $\varepsilon'^{t}\varepsilon=0,1$ respectively.
This motivates the following definition:

\begin{defn}
The integral characteristics $\left[\begin{array}{c}
  \varepsilon \\
  \varepsilon' \\
\end{array}\right]
$ is called \textbf{even} (\textbf{odd}) if $\varepsilon'^{t}\varepsilon=0,1$
\end{defn}

Let us cite the duplication formula for higher dimensional theta functions that will be important in the sequel:
$$
        \Theta^2\left[\begin{array}{c}
               \epsilon \\
               \epsilon_1
             \end{array}\right](0,\tau)=\sum_{\alpha' \in {{\mathbb{Z}_2^g}}}\Theta\left[\begin{array}{c}
                           \epsilon+\alpha' \\
                           2\epsilon_1
                         \end{array}\right](0,2\tau)\Theta\left[\begin{array}{c}
                           \alpha' \\
                           0
                         \end{array}\right](0,2\tau)
$$
Here $\epsilon, \epsilon_1 \in \mathbb{Q}^g$ are any $g$ rational
characteristics. The proof is \cite{Mu} or \cite{RF}. In particular replace $\tau$ with $2\tau$ on both sides of the equation to obtain that:

\begin{equation}\label{eq4}
\Theta^2\left[\begin{array}{c}
               \epsilon \\
               \epsilon_1
             \end{array}\right](0,2\tau)=\sum_{\alpha' \in {{\mathbb{Z}_2^g}}}\Theta\left[\begin{array}{c}
                           \epsilon+\alpha' \\
                           2\epsilon_1
                         \end{array}\right](0,4\tau)\Theta\left[\begin{array}{c}
                           \alpha' \\
                           0
                         \end{array}\right](0,4\tau)
\end{equation}

Let us show the following lemma:
\begin{lem}
Let $\delta \in \mathbb{Z}_2^g$ then:
$$
\Theta\left[\begin{array}{c}
                           0\\
                           \delta \\
                           \end{array}\right](0,\tau)=
\sum_{\beta \in \mathbb{Z}_2^g}exp\left(\pi i \delta \cdot \beta^t\right)\Theta\left[\begin{array}{c}
                          \beta \\
                           0
                         \end{array}\right](0,4\tau)                            $$
\end{lem}

\begin{proof}
We write by the definition of theta functions:
$$
\Theta\left[\begin{array}{c}
                           0\\
                           \delta
                         \end{array}\right](0,\tau)=
                         \sum_{l\in \mathbb{Z}^{g}} exp 2\pi i
\left\{\frac{1}{2}l\tau\l^t+
l^t\frac{\delta}{2}\right\}.
$$
Rewrite the right hand side of the last equation  as :
$$
    \sum_{m \in {Z}^{g},\beta \in Z_2^g} exp 2\pi i
\left\{\frac{1}{2}\left(2m+\beta\right)\tau\left(2m+\beta\right)^t+
(2m+\beta)^t\frac{\delta}{2}\right\}
$$
  Since $exp\left(2\pi i m\delta\right) = 1$ this equals to:
$$
    \sum_{m \in \mathbb{Z}^{g},\beta \in Z_2^g}
   \exp(\pi i \delta^t \beta)\times \exp 2\pi i \left\{\frac{1}{2}\left(m+\frac{\beta}{2}\right)4\tau\left(m+\frac{\beta}{2}\right)^t
\right\}
$$
the last sum equals to :
$$\sum_{\beta \in \mathbb{Z}_2^g}exp\left(\pi i \delta \cdot \beta^t\right)\Theta\left[\begin{array}{c}
                          \beta \\
                           0
                         \end{array}\right](0,4\tau)
$$
by the definition of theta functions.

\end{proof}

\begin{cor} The following identity holds
\[\Theta\left[\begin{array}{c}
                          \beta \\
                           0
                         \end{array}\right](0,4\tau)=
        \frac{1}{2^g}\sum_{\delta\in Z_2^g}\exp\left(-\pi i \delta^t\cdot \beta\right)        \Theta\left[\begin{array}{c}
                           0\\
                           \delta
                         \end{array}\right](0,\tau)
\]
\end{cor}

\begin{proof}
We can treat $\Theta\left[\begin{array}{c}
                          \beta \\
                           0
                         \end{array}\right](0,4\tau)$
as unknowns in a system of linear equations in which the elements of the $2^g\times 2^g$ matrix $A$ are $exp(\pi \delta^t \beta.)$ Multiply $A$ by its Hermitian transpose $A^*.$ Then $(AA^{*})_{ii}=2^g$ since the diagonal element equals: $\sum_\delta \exp(\pi i \delta^t\beta)\exp(-\pi i \delta^t\beta)=2^g.$ if $i\neq j$ then $(AA^*)_{ij}=0.$ This is because this element can be regarded as a sum  of a nontrivial character on the group $Z_2^g.$

\end{proof}

We apply the formula to Eq.~\eqref{eq4}. If $\epsilon_1$ is an integer we  rewrite the equation as:

\begin{equation}
\Theta^2\left[\begin{array}{c}
               \epsilon \\
               \epsilon_1
             \end{array}\right](0,2\tau)=\sum_{\alpha' \in {{\mathbb{Z}_2^g}}}exp\left(\pi i (\epsilon+\alpha')\epsilon_1\right)\Theta\left[\begin{array}{c}
                           \epsilon+\alpha' \\
                           0
                         \end{array}\right](0,4\tau)\Theta\left[\begin{array}{c}
                           \alpha' \\
                           0
                         \end{array}\right](0,4\tau)
\end{equation}

Applying Corollary \textbf{2.4} the last equation equals to :
\begin{equation}
    \sum_{\alpha,\gamma,\delta \in \mathbb{Z}_2^g}\frac{1}{2^{2g}}c_{\epsilon,\epsilon', \alpha,\gamma,\delta}\Theta\left[\begin{array}{c}
               0 \\
               \gamma
             \end{array}\right](0,\tau)\Theta\left[\begin{array}{c}
               0 \\
               \delta
             \end{array}\right](0,\tau)
\end{equation}
and $$c_{\epsilon,\epsilon', \alpha',\gamma,\delta}=\exp\left(\pi i (\epsilon+\alpha')^{t}\epsilon_1\right)\times \exp\left(\pi i(\epsilon'+\alpha')^{t}\gamma\right)\times\exp\left(\pi i \alpha'^{t}\delta\right)
$$
Note that for fixed $\delta,\gamma$ the total coefficient of the product $\Theta\left[\begin{array}{c}
               0 \\
               \gamma
             \end{array}\right](0,\tau)\Theta\left[\begin{array}{c}
               0 \\
               \delta
             \end{array}\right](0,\tau)$ is:
 $$\sum_\alpha\exp\left(\pi i (\epsilon+\alpha')^{t}\epsilon_1\right)\times \exp\left(\pi i(\epsilon'+\alpha')^{t}\gamma\right)\times\exp\left(\pi i \alpha'^{t}\delta\right)
$$
which equals to $0$ unless $\epsilon_1+\gamma+\delta=0.$ In the latter case the coefficient is : $2^g\exp(\pi i \epsilon^t \delta).$ Summarizing we obtain the following corollary:
\begin{cor}
\begin{equation}\Theta^2\left[\begin{array}{c}
               \epsilon \\
               \epsilon_1
             \end{array}\right]^2(0,2\tau)=\sum_{\delta \in Z_2^g}\exp\left(\pi i \epsilon^t \delta \right)\Theta\left[\begin{array}{c}
               0 \\
               \delta             \end{array}\right](0,\tau)\Theta\left[\begin{array}{c}
               0 \\
               \epsilon_1-\delta
             \end{array}\right](0,\tau)
\end{equation}
\end{cor}
We apply the last corollary to explain a procedure to obtain a higher dimensional modular equations
analogues for any prime $p.$

\begin{defn}
Let $D_n$ be a set of vectors $\left[\begin{array}{c}
                                \epsilon \\
                                 \epsilon_1
                                 \end{array}\right]$
such that \begin{itemize}
            \item $\epsilon \in \mathbb{Z}^g, \epsilon_i={0,1}.$
            \item $\epsilon_1 \in \mathbb{Q}^g,\epsilon_{1i}=\frac{l}{2^n}, 0\leq l<2^n$
          \end{itemize}
\end{defn}

For each $\tau \in \mathbf{H_g}$ let
$\alpha_n(\tau)=\Theta\left[\begin{array}{c}
               \epsilon \\
               \epsilon_1
             \end{array}\right](0,\tau)$
and $\left[\begin{array}{c}
               \epsilon \\
               \epsilon_1
             \end{array}\right] \in D_n.$
 $\tau\mapsto \alpha_n(\tau)$ induces a map from $\psi_n(\tau): \mathbf{H_g} \mapsto CP^l$ ($l$ - number of vectors in $D_n$.)
Let $X_n=\psi_n\left(\mathbf{H_g}\right).$  Then $X_1$ is the image in $CP^{2^{2g}-1}$ where $\epsilon , \epsilon_1$ are
integral characteristics. There is a map $\phi_n:X_n
\mapsto X_1$ which omits the non integer
characteristics in the definition of $\psi_n(\tau).$
\begin{lem}
The map $\phi_n$ is a finite map from $X_n \mapsto X_1.$
\end{lem}

\begin{proof}
Because of the Transformation formula for theta functions [RF] under the action of $Sp\left(g,\mathbf{Z}\right)$ , there exists a subgroup $\Delta_n$ of finite index in $Sp\left(g,\mathbf{Z}\right)$ such that $\psi_n(\tau)$ induces a map $\beta_n(\tau)$, $\beta_n:
\mathbf{H_g}/\Delta_n \mapsto X_n$ ( Note: \textbf{$\Delta_n$ is not a congruence subgroup of level $n.$}) Hence, the map $\phi_n$ factors through a map
\[{\bar{\phi_n}}:\mathbf{H_g}/\Delta_n \mapsto \mathbf{H_g}/\Delta_1,\]
which is clearly finite since $\Delta_n$ has a finite index inside $Sp\left(g,\mathbf{Z}\right).$

\end{proof}

Before stating the theorem that describes the map $\phi_n$ more explicitly we state the following definition:
\begin{defn}
Let $H$ be a complex analytic domain. $f_1...f_n:H\mapsto C$ be complex analytic functions. We call $f$ constructible from $f_1...f_n$ if
\begin{itemize}
   \item $f$ is algebraic above $\mathbb{C}(f_1...f_n)$
   \item The Galois group of $C(f)$ above $f_1...f_n$ is solvable.
   equivalently $f$ can be expressed through radical expressions involving $f_1...f_n.$
 \end{itemize}
\end{defn}

\begin{thm}
Let $\epsilon_1 \in \mathbb{Q}^g,\epsilon_{1i}=\frac{l}{2^n}.$ Let
$\epsilon \in Z_2^g.$ Then $\Theta\left[\begin{array}{c}
               \epsilon \\
               \epsilon_1
             \end{array}\right](0,\tau)$
 is constructible  from $\Theta\left[\begin{array}{c}
               \eta\\
               \eta_1
             \end{array}\right](0,\tau)$
and $\eta, \eta_1$ are integral characteristics.
\end{thm}


\begin{proof}
Assume inductively that the theorem is true for all
characteristics $\Theta\left[\begin{array}{c}
               \delta \\
               \delta_1
             \end{array}\right](0,\tau)$ such that $\delta \in
\mathbb{Z}^g$ and $\delta_{1i}=\frac{l}{2^{n-1}}.$ The duplication
formula implies:
$$
        \Theta^2\left[\begin{array}{c}
               \epsilon \\
               \epsilon_1
             \end{array}\right](0,\tau)=\sum_{\alpha' \in {{\mathbb{Z}_2^g}}}\Theta\left[\begin{array}{c}
                           \epsilon+\alpha' \\
                           2\epsilon_1
                         \end{array}\right](0,2\tau)\Theta\left[\begin{array}{c}
                           \alpha' \\
                           0
                         \end{array}\right](0,2\tau)
$$
But
$$\Theta\left[\begin{array}{c}
                           \epsilon+\alpha' \\
                           2\epsilon_1
\end{array}\right](0,2\tau)$$
satisfies the induction hypothesis. So it is constructible from
$$
        \Theta^2\left[\begin{array}{c}
               \eta \\
               \eta_1
             \end{array}\right](0,2\tau)$$ Apply the formulas from corollary \textbf{2.6} to see that $$\Theta^2\left[\begin{array}{c}
               \eta \\
               \eta_1
             \end{array}\right](0,2\tau)$$ is constructible from $$\Theta\left[\begin{array}{c}
               \eta \\
               \eta_1
             \end{array}\right](0,\tau).$$
 Hence $$\Theta\left[\begin{array}{c}
               \epsilon \\
               \epsilon_1
             \end{array}\right](0,\tau)$$ is constructible.
\end{proof}


Note that the proof gives slightly more i.e. a recursive process how to construct  the expressions $$\Theta\left[\begin{array}{c}
               \epsilon \\
               \epsilon_1
             \end{array}\right](0,\tau)$$ from
             $$\Theta\left[\begin{array}{c}
               \eta \\
               \eta_1
             \end{array}\right](0,\tau)$$
We reformulate it in the following corollary:
\begin{cor}
Let $P \in X_1$ and let $Q \in \phi_n^{-1}(P)$ then there exists a process
that produces an algebraic relationship between the
coordinates of $Q$ and coordinates of $P.$
\end{cor}

We rely on the last theorem and the methods developed in [Ko] to obtain process to produce generalized modular equation in  any dimension.
\begin{defn}
Let $p$ be a prime. A modular equation of order $p$ will be any non trivial polynomial identity between $\Theta\left[\begin{array}{c}
               \eta \\
               \eta_1
             \end{array}\right](0,\tau)$ and  $\Theta\left[\begin{array}{c}
               \eta \\
               \eta_1
             \end{array}\right](0,p\tau)$
$\eta,\eta_1$ integral characteristics.
\end{defn}

To introduce the theorem from [Ko] we need the notion of a function order $k$ of characteristics  and characteristics $\left[\begin{array}{c}
         0...0 \\
         0...0
 \end{array}\right]$.
\begin{defn}
$f:C^g\times\mathbf{H_g}\mapsto C$ is an analytic function of order $k$ and characteristics $\left[\begin{array}{c}
         0...0 \\
         0...0
 \end{array}\right]$ if the following relation is satisfied:

$$f\left(z + n + m^t\tau, \tau\right) =  \exp\left\{2\pi i(-km^tz-km^t\tau m)\right\}f\left(z,\tau\right)$$
\end{defn}

Let $k=p_1p_2$ and
$p_1,p_2$ are even arbitrary numbers. We quote the following theorem
from \cite{Ko}.

\begin{thm}
Let $f$ be a function of characteristics
 $\left[\begin{array}{c}
         0...0 \\
         0...0
 \end{array}\right]$ and even order $k.$ If $\left[\begin{array}{c}
                           \mu \\
                           \mu'
                         \end{array}\right]$ is an integral odd
                         characteristics then the following
                         identity is valid :

$$
\sum_{\nu,\nu', 0 \leq \nu_i\leq p_1, 0 \leq \nu'_i \leq
p_2}(-1)^{\mu \nu- \mu' \nu'}f\left(\frac{\nu}{p_1}+\tau
\frac{\nu'}{p_2}\right)=0
$$

and $\mu \nu=\sum_i \mu_i \nu_i.$
\end{thm}
To obtain a modular equation for an odd number $p,$ Choose
$k=2^{[log_2(p)]+1} , l=k-p$ and examine the function :
$$
    f=\Theta\left[\begin{array}{c}
                           0 0...0 \\
                           0 0...0
                         \end{array}\right]^l(z,\tau)\Theta\left[\begin{array}{c}
                           0 0...0 \\
                           0 0...0
                         \end{array}\right](pz,p\tau)
$$
This function is of characteristics  $\left[\begin{array}{c}
                           0 0...0 \\
                           0 0...0
                         \end{array}\right]$ and order $k.$
Define $p_1 = 2^{[log_2(p)]}, p_2=2.$ Apply Theorem 2.14 to obtain the following:
\begin{thm} (\it Modular equation for $p$)
For any $\left[\begin{array}{c}
                           \mu \\
                           \mu'
                         \end{array}\right]$ odd integral characteristics:
$$
    \sum_{\nu,\nu', 0 \leq \nu'_i\leq p_1, 0 \leq \nu_i \leq 1
}(-1)^{\mu \nu- \mu' \nu'} \Theta\left[\begin{array}{c}
                           \nu_i \\
                           \frac{2\nu'_i}{p_1}
                         \end{array}\right]^l(0,\tau)\Theta\left[\begin{array}{c}
                           \nu_i \\
                           \frac{2\nu'_i}{p_1}
                         \end{array}\right](0,p\tau)=0
$$
\end{thm}
where $p_1=2^{[log_2(p)]}.$ The main theorem of the paper is given below:

\begin{thm}
There exist explicit equation connecting
\[
\Theta\left[\begin{array}{c}
                           \chi_i \\
                           \chi'_i
\end{array}\right]^l(0,\tau), \textit{ and }
\Theta\left[\begin{array}{c}
                           \chi_i \\
                           \chi'_i
\end{array}\right]^l(0,p\tau)
\]
$\chi_i, \chi'_i$ are $g$ integral characteristics.
\end{thm}


\begin{proof}
According to theorem 2.15:$$
    \sum_{\nu,\nu', 0 \leq \nu'_i\leq p_1, 0 \leq \nu_i \leq 1
}(-1)^{\mu \nu- \mu' \nu'} \Theta\left[\begin{array}{c}
                           \nu_i \\
                           \frac{2\nu'_i}{p_1}
                         \end{array}\right]^l(0,\tau)\Theta\left[\begin{array}{c}
                           \nu_i \\
                           \frac{2\nu'_i}{p_1}
                         \end{array}\right](0,p\tau)=0$$ and  $\left[\begin{array}{c}
                           \mu \\
                           \mu'
                         \end{array}\right]$  is an odd characteristics.
  Applying  theorem 2.9 conclude $\Theta\left[\begin{array}{c}
                           \nu_i \\
                           \frac{2\nu'_i}{p_1}
                         \end{array}\right]^l(0,\tau)$ is constructible from $\Theta\left[\begin{array}{c}
               \eta\\
               \eta_1
             \end{array}\right](0,\tau)$ and $\eta,\eta_1$ is an integral characteristics. Replace $\Theta\left[\begin{array}{c}
                           \nu_i \\
                           \frac{2\nu'_i}{p_1}
                         \end{array}\right]^l(0,\tau), \Theta\left[\begin{array}{c}
                           \nu_i \\
                           \frac{2\nu'_i}{p_1}
                         \end{array}\right]^l(0,p\tau)$
 with the corresponding radical expression involving $\Theta\left[\begin{array}{c}
               \eta\\
               \eta_1
\end{array}\right](0,\tau)$  to conclude the result.

\end{proof}


\section{The one dimensional case}
In this section we explain the connection between the theory developed in the last section and the usual one dimensional theory of modular equations.

Let $E$ be an elliptic curve given in Legendre's normal form $y^2=x(x-1)(x-\lambda).$ if $\tau$ denotes the period that is induced by $E,$ we have the following expression for $\lambda$ as function of $\tau:$

\begin{equation}
    \lambda(\tau) = \frac{\Theta^4\left[\begin{array}{c}
               1 \\
               0
             \end{array}\right](0,\tau)}{\Theta^4\left[\begin{array}{c}
               0 \\
               0
             \end{array}\right](0,\tau)}
\end{equation}
Recall the identity

\begin{equation}
\Theta^4\left[\begin{array}{c}
               0 \\
               0
             \end{array}\right](0,\tau)=\Theta^4\left[\begin{array}{c}
               1 \\
               0
             \end{array}\right](0,\tau) + \Theta^4\left[\begin{array}{c}
               0 \\
               1
             \end{array}\right](0,\tau)
\end{equation}
dividing both sides by $\Theta^4\left[\begin{array}{c}
               0 \\
               0
             \end{array}\right](0,\tau)$ we see that

\begin{equation}
    1-\lambda(\tau) = \frac{\Theta^4\left[\begin{array}{c}
               0 \\
               1
             \end{array}\right](0,\tau)}{\Theta^4\left[\begin{array}{c}
               0 \\
               0
             \end{array}\right](0,\tau)}
\end{equation}
Now set $\Theta^4\left[\begin{array}{c}
               0 \\
               0
             \end{array}\right](0,\tau)=\theta_0(\tau), \Theta^4\left[\begin{array}{c}
               0 \\
               1
             \end{array}\right](0,\tau)=\theta_1(\tau),\Theta^4\left[\begin{array}{c}
               1 \\
               0
             \end{array}\right](0,\tau)=\theta_2(\tau)$
We write:
\begin{equation}
    \sqrt[4]{\lambda(\tau)}=\frac{\theta_1(\tau)}{\theta_0(\tau)},
\end{equation}
  and
\begin{equation}
 \sqrt[4]{1-\lambda(\tau)}=\frac{\theta_2(\tau)}{\theta_0(\tau)}.
 \end{equation}
The proof of the main theorem in the last section applied to the one
dimensional case produces a homogenous radical expression of the form:
\begin{equation}
F\left(\theta_i(\tau)\theta_j(p\tau)\right)=0.
\end{equation}
Divide each term of the expression by $\theta_0(\tau)\theta_0(p\tau).$ Using the definition of   $\lambda(\tau)$ replace its quotient by $\lambda(\tau)$ and $\lambda(p\tau)$ respectively. We obtain the following classical theorem:

\begin{thm}
There exists an algebraic relation between $\lambda(\tau)$ and
$\lambda(p\tau).$
\end{thm}
Note that the proof of this theorem we obtain does not use the usual modular group theory. The proof is constructive and provides an alternative way to construct modular equations for any $p$.

As an example consider the case $p=3$ then applying the algorithm we
obtain:
\begin{equation}
    \theta_0(\tau)\theta_0(3\tau)=\theta_1(\tau)\theta_1(3\tau)+\theta_2(\tau)\theta_2(3\tau)
\end{equation}
divide the two sides of the last equation by $\theta_0(\tau)\theta(3\tau)$ we obtain the classical modular
equation:
\[1=\sqrt[4]{\lambda(\tau)\lambda(3\tau)}+\sqrt[4]{\left(1-\lambda(\tau)\right)\left(1-\lambda(3\tau)\right)}.\]
More details can be found in [Bo].

\section{Examples}

We apply our theory to two cases as an example:

\subsection{p=3, g=2}
In this section we outline the result of the method applied above to
$p=3$ for $g=2.$  First we observe that for $g=2$ we have 6 odd and
10 even characteristics. We conclude immediately that overall we
have 6 modular equations for each $p.$ Using the Theorem 2.7 we
write for $p=3$ :
\begin{thm}
    For any $\left[\begin{array}{c}
                           \mu \\
                           \mu'
                         \end{array}\right]$ odd integral characteristics the
                         following identities are true:
$$
    \sum_{\nu,\nu', 0 \leq \nu'_i\leq 1, 0 \leq \nu_i \leq 1
}(-1)^{\mu \nu- \mu' \nu'} \Theta\left[\begin{array}{c}
                           \nu_i \\
                           \nu'_i
                         \end{array}\right]^l(0,\tau)\Theta\left[\begin{array}{c}
                           \nu_i \\
                           \nu'_i
                         \end{array}\right](0,3\tau)=0
$$

\end{thm}
To achieve a more compact set of identities we rely on the
classification of identities of power 4 achieved in \cite{AK}.
\begin{defn}
Let $$0 =\left[\begin{array}{c}
  0  \\
  0  \\
\end{array}\right], 1 = \left[\begin{array}{c}
  1 \\
  0 \\
\end{array}\right], 2 = \left[\begin{array}{c}
  0 \\
  1 \\
\end{array}\right], 3=\left[\begin{array}{c}
  1 \\
  1 \\
\end{array}\right].$$

\end{defn}
Using the last definition we can write any 2 dimensional
characteristics using the vectors above. For example:
$$
03 =\left[\begin{array}{c}
  0 1  \\
  0 1  \\
  \end{array}\right]
$$
We define matrix $A$ the encodes even characteristics in the
following way :
$$
    A= \left(
\begin{array}{ccc}
  11 & 01 & 10 \\
  22 & 20 & 02 \\
  33 & 21 & 22 \\
\end{array}\right),
\left(%
\begin{array}{c}
  00 \\
\end{array}%
\right)
$$
Then the following classification of 2 dimensional theta identities
is given in \cite{AK}

\begin{thm}
There are two types of relations for theta functions in power 4:
\begin{itemize}
    \item Type I - corresponds to generalized diagonals of matrix
    A on the one side and the characteristics $\left(%
\begin{array}{c}
  00 \\
\end{array}%
\right)$ on the other for example the following identity is true :
$$\Theta^4_{00}= \Theta^4_{11}+\Theta^4_{20}+\Theta^4_{12}(B)
$$
    \item Type II - Correspond to $2 \times 2$ sub matrices of $A$
    putting 2 not the same row and not the same column entries on
    each side of the equation. For example:
    $$ \Theta^4_{11} + \Theta^4_{20}=\Theta^4_{01} + \Theta^4_{22}$$
\end{itemize}
\end{thm}
The identities of power 4 of theta constants were obtained applying
the Gauss elimination procedure to the matrix that defines the power
4 identities of theta functions \cite{AK}. Since the same matrix
governs the identities involving terms of the form
$\Theta(0,\tau)\Theta(0, 3\tau)$ an immediate corollary of the last
theorem is the following classification of identities involving
prime $p=3$:
\begin{thm}
There are two types of relations for theta functions involving the prime
3:
\begin{itemize}

    \item Type I - corresponds to generalized diagonals of matrix
    A on the one side and the characteristics $\left(%
\begin{array}{c}
  00 \\
\end{array}%
\right)$  for example the following identity is true :
\begin{small}
\[
\Theta_{00}(0, \tau )\Theta_{00}(0, 3\tau ) = \Theta_{11}(0, \tau)\Theta_{11}(0, 3\tau)
+\Theta_{20}(0,\tau)\Theta_{20}(0,3\tau)+\Theta_{12}(0,\tau)\Theta_{12}(0,3\tau)
\]
\end{small}
    \item Type II - Correspond to $2 \times 2$ sub matrices of $A$
    putting 2 not the same row and not the same column entries on
    each side of the equation. For example:
\begin{small}
\[
\Theta_{11}(0,\tau)\Theta_{11}(0,3\tau)  + \Theta_{20}(0,\tau)\Theta_{20}(0,3\tau)
    =\Theta_{01}(0,\tau)\Theta_{01}(0,3\tau) + \Theta_{22}(0,\tau)\Theta_{22}(0,3\tau)
\]
\end{small}

\end{itemize}
\end{thm}

This gives an algebraic way to evaluate
$\Theta_{ij}(0,3\tau)$ from  $\Theta_{ij}(0,\tau).$ Compare with
\cite{CKL}.

\subsection{p=7, g=2}Let us write the equations explicitly for
$p=7.$ According to the recipe outlined above we need to choose
$k=8$ and thus our basic function will be:
$$
    f=\Theta\left[\begin{array}{c}
                           0 0 \\
                           0 0
                         \end{array}\right](z,\tau)\Theta\left[\begin{array}{c}
                           0 0 \\
                           0 0
                         \end{array}\right](7z,7\tau)
$$
for $p=7$  Now we use theorem 2.14 to obtain
the 6 equations for each odd characteristics. For example if the  odd
characteristics is $$\left[\begin{array}{c}
                           1 0 \\
                           1 0
                         \end{array}\right]$$
The equation is:
$$
\sum_{0 \leq \nu_1,\nu_2 \leq 1, 0 \leq \nu_1', \nu_2' \leq
2}\left(-1\right)^{\nu_1-\nu_1'}\Theta \left[\begin{array}{c}
                           \nu_1  \nu_2 \\
                           \frac{\nu_1'}{2}  \frac{\nu_2'}{2}
                            \end{array}\right](0,\tau)\Theta\left[\begin{array}{c}
                           \nu_1  \nu_2 \\
                           \frac{\nu_1'}{2}  \frac{\nu_2'}{2}
                         \end{array}\right](0,7\tau)=0
$$

To translate the equation into equation involving integral
characteristics we use the duplication formulas. For example we
write:
$$
    \Theta^2\left[\begin{array}{c}
                           1  1 \\
                           \frac{1}{2}  \frac{1}{2}
                            \end{array}\right](0,\tau)=\Theta\left[\begin{array}{c}
                            1 1\\
                            1 1\end{array}\right](0,2\tau)\Theta\left[\begin{array}{c}
                            0 0\\
                            0 0\end{array}\right](0,\tau)+\Theta\left[\begin{array}{c}
                            0 0\\
                            1 1\end{array}\right](0,2\tau)\Theta\left[\begin{array}{c}
                            1 1\\
                            0 0\end{array}\right](0,\tau)$$

Since the other two theta functions with integral characteristics are equal to 0. Similar formulas hold for the other functions of the form : $\Theta\left[\begin{array}{c}
                           \nu_1  \nu_2 \\
                           \frac{\nu_1'}{2}  \frac{\nu_2'}{2}
                            \end{array}\right](0,\tau)$
Substituting we obtain formulas involving integral characteristics at point $\tau, 2\tau, 14\tau.$ To reduce to equations involving $\tau, 7\tau$ we apply the formulas from corollary 2.5 to
functions
$$\Theta\left[\begin{array}{c}
                           \nu_1  \nu_2 \\
                           \nu_1'  \nu_2'
                            \end{array}\right](0,2\tau)
$$ and $\nu_i, \nu_i'$ are integral characteristics.

\vspace{1cm}
\noindent Yaacov Kopeliovich \\
540 Madison Avenue, 6 -th floor\\
New York NY 10022\\
Email: ykopeliovich@yahoo.com,

\hspace{.7cm} ykopeliovich@medtolife.com


\begin{thebibliography}{12}
    \bibitem [AK]{AK} R.Adin, Y.Kopeliovich,  Short Eigenvectors and
    Multidimensional Theta Functions, {\em Linear Algebra and
    Appl.} \textbf{257}(1)(1997)  49-63

    \bibitem [Bo]{Bo} J. Borwein, P.Borwein, {\em $\Pi$ and the AGM},
     A Wiley Interscience publication, 1987

    \bibitem[CKL]{CKL} R.Carls, D.Kohel and D.Lubicz,  Higher Dimensional
    3-Adic CM Construction {\em Preprint}

    \bibitem[CL]{CL} R.Carls, and D.Lubicz,  A $p$-adic quasi quadratic point counting algorthm {\em \url{http://arxiv.org/PS_cache/arxiv/pdf/0706/0706.0234v2.pdf}}

    \bibitem [FK1]{FK1} H. Farkas, Y. Kopeliovich, New Theta Constant
    Identities {\em Israel Journal of Mathematics }\textbf{82}(1)(1993) 133-140

    \bibitem [FK2]{FK2} H. Farkas, Y.Kopeliovich,  New Theta Constant
    Identities II {\em Proceeding of AMS.}\textbf{123}(4)(1995) 1009-1020

    \bibitem[GS]{GS} P.Gaudry, E.Schost, Modular Equations for Hyperelliptic curves {\em \url{http://www.csd.uwo.ca/~eschost/publications/papier2.pdf}}

    \bibitem [Ko]{Ko} Y. Kopeliovich,  Multi Dimensional Theta Constant
    Identities {\em Journal of Geometric Analysis} \textbf{8} (4)(1998) 571-581

    \bibitem[Ma]{Ma} M.Madsen, A general framework for $p$ - adic
    point counting and applications to elliptic curves on Legendre
    form {\em
    \url{http://www.imf.au.dk/publications/pp/2004/imf-pp-2004-2.pdf}}

    \bibitem[Me]{Me} J.-F.Mestre, Notes on Talk given at seminar of Cryptography
    at Rennes 2002. {\em
    \url{http://www.math.univ-rennes1.fr/crypto/2001-02/mestre.ps}}

    \bibitem [Mu]{Mu} D. Mumford, {\em Tata Lectures on Theta II} (Progress in
    Mathematics, Birkhauser 1984)

    \bibitem[RF]{RF} H. Rauch and H.Farkas {\em Theta functions with application to Riemann Surfaces} (William and Wilkins
    Balt. Md. 1974)

\end{thebibliography}
\end{document}